\documentclass[%
  twocolumn
   , colorlinks 
   , hidempi
   , 
]{mpi2015-cscpreprint}

\usepackage[american]{babel}

\usepackage{graphicx}
\usepackage{todonotes}
\usepackage{orcidlink}
\usepackage{rorlink}

\usepackage{amssymb}
\usepackage{amsthm}

\AtBeginDocument{

}

\usepackage[vlined,longend,linesnumbered,ruled]{algorithm2e}
\newcommand{\tran}{\ensuremath{\mkern-1.5mu\mathsf{T}}}

\usepackage{booktabs}
\usepackage{multirow}

\usepackage{pgfplotstable}
\usepackage{pgfplots}
  \pgfplotsset{compat = 1.18,
    unbounded coords = jump,
    every axis plot/.append style={line width=.75pt}
  }
  \usepgfplotslibrary{external}
  \tikzset{external/system call = {%
    pdflatex \tikzexternalcheckshellescape%
      -halt-on-error
      -interaction=batchmode
      -jobname "\image" "\texsource"}}
  \tikzexternaldisable%

\definecolor{mycolor1}{HTML}{A6CEE3}
\definecolor{mycolor2}{HTML}{1F78B4}
\definecolor{mycolor3}{HTML}{B2DF8A}
\definecolor{mycolor4}{HTML}{33A02C}

\pgfplotscreateplotcyclelist{mylist}{%
{mycolor1, dashed},
{mycolor2, dash pattern = on 2pt off 2pt on 1pt off 2pt on 1pt off 2pt},
{mycolor3, densely dashed},
{mycolor4, dotted},
}
\pgfplotsset{
    every axis/.append style={
        cycle list name=mylist
    }
}
\newlength{\figwidth}
\newlength{\figheight}

\usepackage{caption}
\usepackage{subcaption}

\begin{document}


\title{On the Solution of Large-scale Non-autonomous Differential Riccati Equations: a Numerical Study}

\author[$\ast$]{Eugenio Mancinelli\orcidlink{0009-0006-8022-6940}}
\affil[$\ast$]{Dipartimento di Matematica,
Alma Mater Studiorum Universit\`a di Bologna~\rorlink{https://ror.org/01111rn36},
  \email{eugenio.mancinelli3@unibo.it}, }

\author[$\dagger$]{Davide Palitta\orcidlink{0000-0002-6987-4430}}
\affil[$\dagger$]{Dipartimento di Matematica and (AM)$^2$,
Alma Mater Studiorum Universit\`a di Bologna~\rorlink{https://ror.org/01111rn36},
  \email{davide.palitta@unibo.it}}

\author[$\ddagger$]{Jens Saak\orcidlink{0000-0001-5567-9637}}
\affil[$\ddagger$]{Research Group Computational Methods in Systems and Control Theory (CSC), Max Planck
Institute for Dynamics of Complex Technical Systems~\rorlink{https://ror.org/030h7k016},
  \email{saak@mpi-magdeburg.mpg.de}}

\shorttitle{Solution of DREs: a Numerical Study}
\shortauthor{E. Mancinelli, D. Palitta, J. Saak}
\shortdate{}



\abstract{We explore the numerical solution of large-scale non-autonomous Differential Riccati Equations (DREs). While we assume to discretize the differential operator using a Backward Differentiation Formula (BDF) of order $s$, we solve the generalized Algebraic Riccati Equation (gARE) resulting at each time step by different state-of-the-art methods. In particular, we compare the performance of the inexact Newton-Kleinman method with line search and the low-rank RADI iteration, considering for both methods two different initialization strategies: zero initialization and warm-start. A comprehensive panel of numerical results illustrate the potential and limitations of these methods when employed within a numerical pipeline for the solution of DREs, rather than for the isolated solution of a single gARE, as commonly considered in the existing literature.
}

 \novelty{This paper provides a fresh numerical comparison about the performance achieved by different state-of-the-art solvers for gAREs when these are adopted as inner solver for the solution of non-autonomous DREs. While available work surveys the performance of these methods in the solution of a single gARE, facing long sequences of gAREs as the ones stemming from the discretization of DREs by BDFs unveils new drawbacks and potential.}

\maketitle


\section{Introduction}%
\label{sec:intro}

The goal of this paper is to provide a numerical study comparing different state-of-the-art solvers for generalized algebraic Riccati equations arising in the step equations of Backward Differentiation Formula (BDF) methods for Differential Riccati Equations (DREs)
{\small
\begin{equation}\label{eq:main}
    \begin{array}{rl}
         -E^{\tran}\dot X(t)E =&\hspace{-8pt} E^{\tran}X(t)\left(A+\dot E\right) + \left(A+\dot E\right)^{\tran}X(t)E\\
         &+ C^{\tran}QC - E^{\tran}X(t)BR^{-1}B^{\tran}X(t)E, \\
X({t_f}) =&\hspace{-8pt} X_{t_f},\\
    \end{array}
\end{equation}}
where $t\in[t_0,t_f]$ represents time on a given interval. The coefficients $A,E\in\mathbb{R}^{n\times n}$, $B\in\mathbb{R}^{n\times p}$, $C\in\mathbb{R}^{q\times n}$ are full-rank matrices with $p+q\ll n$. Moreover, $Q\in\mathbb{R}^{q\times q}$ and $R\in\mathbb{R}^{p\times p}$ are symmetric positive-definite\footnote{For $Q$ actually semi-definiteness is sufficient.}.

We consider both the cases of autonomous and non-autonomous equations. In the latter case, not only the solution $X$ but also the coefficient matrices $A$, $E$, $B$, and $C$ depend on time, although this dependence is not explicitly indicated in the formulation~\eqref{eq:main}. Throughout the paper, however, the matrices $Q$ and $R$ are assumed to be constant in time.

The most prominent application of this equation, and also the main reason for formulating it backward in time, is the linear-quadratic optimal control of linear time-varying (in the non-autonomous case) or linear time-invariant (in the autonomous case) dynamical systems on a finite time horizon, i.e. $t_f<\infty$. There, using a specific ansatz for the dual state, the DRE replaces the dual equation in the optimization problem, explaining the reversion of time. However, after a change of variable, equation~\eqref{eq:main} can be expressed as an initial value problem; see, e.g.,~\cite[Section 2]{baranetAl24}. In the remainder of this work we assume that such a change of variable has been performed so that we can adopt a forward-in-time approach. In the context of linear-quadratic regulator problems, the matrices $Q$ and $R$ in~\eqref{eq:main} represent the penalty matrices considered in the cost function that needs to be minimized to compute the optimal control. This different nature of $Q$ and $R$ explains why we consider these matrices to be constant in time, in contrast to $A$, $E$, $B$, and $C$, that are the actual coefficient matrices of the underlying linear system and thus possibly time-dependent.


A classic numerical pipeline for the solution of~\eqref{eq:main} sees first the discretization of the differential problem on a given time grid $t_0\leq t_1\leq\cdots\leq t_\ell\equiv t_f$ and then the solution of a sequence of algebraic problems, in a time-marching fashion. For the sake of simplicity and to better identify the effect of the different inner solvers, we restrict ourselves to the case of a fixed time grid with constant step length. The nature of the discrete problems, as well as the number of subproblems to be solved at each time step, depends on the adopted discretization scheme. Several options are available in the literature, ranging from  midpoint and trapezoidal rules~\cite{LANGetAl2015}, to Rosenbrock methods~\cite{BennerMena2013} and Backward Differentiation Formulas~\cite{BenM04}.

In this work, we focus on the latter class of time-integra\-tors. When BDFs are applied to~\eqref{eq:main}, we need to solve at each time step $t_{k}$ a generalized Algebraic Riccati Equation (gARE) of the form
\begin{equation}\label{eq:gAREs}
    \begin{split}
  \widehat C_{k}^{\tran}\widehat Q_{k}\widehat C_{k} + \widehat A_{k}^{\tran}X_{k}\widehat E_{k}+
 \widehat E_{k}^{\tran} X_{k}\widehat A_{k}\\
  -\widehat E_{k}^{\tran}X_{k}\widehat B_{k}\widehat R_{k}^{-1}\widehat B_{k}^{\tran}X_{k}\widehat E_{k}=0,
\end{split}
\end{equation}
where $X_{k}\approx X(t_{k})$, and the expression of the coefficient matrices depends on the order of the selected BDF\@; more details on the construction of these matrices will be given in Section~\ref{sec:BDFs}. A similar study for Rosenbrock schemes requiring the solution of Lyapunov equations for the algebraic subproblems can be found in~\cite{SchS24}.

Assuming we split the time interval $[t_0,t_f]$ in $\ell-1$ subintervals, the solution of the $\ell$ gAREs~\eqref{eq:gAREs} is often the computational bottleneck of the solution of DREs.
Our goal is to test different state-of-the-art low-rank solvers for gAREs and measure their performance in terms of the overall solution process required by~\eqref{eq:main}. Notice that this is a rather different setting than testing solvers for gAREs only\footnote{A thorough numerical comparison of solvers for gAREs can be found in~\cite{BenneretAl2020}.}.
Indeed, when facing non-autonomous DREs, a list of computational challenges can come up. These are due to different factors that may be all present in the DRE of interest. First, the problem dimension $n$ can be large so that dense numerical linear algebra tools are out of the picture. Second, the dynamics encoded in~\eqref{eq:main} may require high resolution to be fully captured by the numerical approximations. This means that a fine time grid is needed and the number of gAREs~\eqref{eq:gAREs} to be solved grows accordingly. Third, due to the nature of the problem at hand or the adopted solver, the inner-factor dimension of the iterates may become large. More precisely, our low-rank solvers provide an approximate solution of the form $L_kD_kL_k^\top\approx X(t_k)$ where $L_k\in\mathbb{R}^{n\times r_k}$ and $D_k\in\mathbb{R}^{r_k\times r_k}$. The magnitude of $r_k$ is an important factor that should not get overlooked. Indeed, dealing with objects with large $r_k$ can potentially jeopardize the effectiveness of any low-rank solver.

The state-of-the-art solvers for gAREs we are going to compare are (i) the inexact Newton-Kleinman method (iNK) with line search~\cite{BenneretAl2016,FeitzingeretAl2009,Kleinman1968,BennerByers1998}, and (ii) the RADI scheme~\cite{BenneretAl2018}.
Other numerical schemes for algebraic Riccati equations can be found in the literature; see, e.g., Krylov methods~\cite{Simoncini2016}. However, their use in case of generalized equations ($\widehat E_{k}\neq I$ in~\eqref{eq:gAREs}) is not well assessed.
Different approaches for DREs, especially in the class of Krylov-projection methods, include a first dimension reduction step applied to~\eqref{eq:main}, so as to obtain a smaller DRE to solve; see, e.g.,~\cite{KirstenSimoncini2020,BehBH21}. However, the potentially non-autonomous nature of the problems we are interested in makes the identification of suitable approximation spaces rather difficult, and this aspect has not been solved in the literature so far.

Both iNK and RADI are iterative methods and thus require a starting guess $X_{k}^{(0)}$.
When dealing with a single gARE, $X^{(0)}$ is often chosen as the zero matrix. In our context, justified by the time evolution nature of our problem, we will also report on the performance of iNK and RADI when $X_{k-1}$, namely the approximate solution computed at time $t_{k-1}$, is adopted as initial guess for the $k$-th gARE\@.

Here is a synopsis of the paper. In Section~\ref{sec:BDFs} we revisit the use of BDFs for DREs to properly define all the matrices involved in the algebraic problem~\eqref{eq:gAREs} at each time step. Section~\ref{sec:solver4gAREs} revises the state-of-the-art solvers for gAREs we test in this work, namely the inexact Newton-Kleinman method with line search (Section~\ref{sec:NK}) and RADI (Section~\ref{sec:RADI}). The main contribution of the paper is in Section~\ref{sec:numerics}, where we report the results of our comprehensive numerical testing. The paper ends with Section~\ref{sec:concl}, where we draw our conclusions.

\section{BDFs for DREs}\label{sec:BDFs}

\textit{Backward Differentiation Formulas} (BDFs)~\cite{Ascher1998} are a family of implicit methods used to numerically solve ordinary differential equations. The idea behind these methods is to express the solution at a given time in terms of the already-computed solutions from previous time steps.

The general BDF of order $s$ allows us to discretize~\eqref{eq:main} as
\begin{equation}\label{eq:BDF}
    X_{k}=\sum_{j=1}^s-\alpha_{j}X_{k-j}+\tau\beta\mathcal{R}(t_{k},X_{k}),
\end{equation}
where $\tau=(t_{f}-t_{0})/(\ell-1)$ is the time step size, whereas
the expression of the scalars $\alpha_j$'s and $\beta$ depends on the selected order $s$. Table~\ref{tab:BDF coeff} reports these values for $s\leq 6$, as BDF schemes become unstable for larger orders; see, e.g.,~\cite{Ascher1998}.
\tabcolsep1.5ex
\begin{table}[h!]
    \centering
    {\renewcommand{\arraystretch}{1.5}
    \begin{tabular}{c|c|c|c|c|c|c|c}
$s$ & $\beta$ & $\alpha_1$ & $\alpha_2$ & $\alpha_3$ & $\alpha_4$ & $\alpha_5$ & $\alpha_6$  \\
\hline
1 & $1$ & $-1$ & & & & &  \\
2 & $\frac{2}{3}$ & $-\frac{4}{3}$ & $\frac{1}{3}$ & & &  & \\
3 & $\frac{6}{11}$ & $-\frac{18}{11}$ & $\frac{9}{11}$ & $-\frac{2}{11}$ & &  & \\
4 & $\frac{12}{24}$ & $-\frac{48}{25}$ & $\frac{36}{25}$ & $-\frac{16}{25}$ & $\frac{3}{25}$ &  & \\
5 & $\frac{60}{137}$ & $-\frac{300}{137}$ & $\frac{300}{137}$ & $-\frac{200}{137}$ & $\frac{75}{137}$ & $-\frac{12}{137}$ &  \\
6 & $\frac{60}{147}$ & $-\frac{360}{147}$ & $\frac{450}{147}$ & $-\frac{400}{147}$ & $\frac{225}{147}$ & $-\frac{72}{147}$ & $\frac{10}{147}$ \\
    \end{tabular}}
    \caption{$\alpha_j$'s and $\beta$ determining the $s$-step BDF~\eqref{eq:BDF}.}%
\label{tab:BDF coeff}
\end{table}
In~\eqref{eq:BDF}, $\mathcal{R}(t_{k},X_{k})$ is a shorthand notation for the right-hand side of the first equation in~\eqref{eq:main}, where the coefficient matrices are evaluated at time $t_{k}$ and $X(t_{k})$ is replaced by its approximation $X_{k}$.

By using an $s$\/th-order BDF and assuming the computed solutions to be in an $LDL^{\tran}$-format, namely 
\[X_{k-j}=L_{k-j}D_{k-j}L_{k-j}^{\tran},\quad \text{for }j=0,\ldots,s-1,\]
a simple manipulation of the terms shows that, at each time step, we are left with solving a gARE~\eqref{eq:gAREs}, where
\begin{eqnarray*}
\label{eq:bdf_discr}
    \widehat C_{k}^{\tran}:=&\left[C(t_{k})^{\tran},E(t_{k})^{\tran}L_k,\ldots,E(t_{k})^{\tran}L_{k-s+1}\right],\\
    \widehat Q_{k}:=&\begin{bmatrix}
        \tau\beta Q &  &  &  \\
 & -\alpha_1 D_k &  &  \\
 &  & \ddots &  \\
 &  &  & -\alpha_s D_{k-s+1}
    \end{bmatrix},\\
        \widehat A_{k}:=&\tau\beta\left(A(t_{k})-\dot E(t_{k})\right)-\frac{1}{2}E(t_{k}),
    \end{eqnarray*}
    with $\widehat E_{k}:=E(t_{k})$, $\widehat B_{k}:=B(t_{k})$, and $\widehat R_{k}^{-1}:=\tau\beta R^{-1}$.
    Observe that to start a BDF of order $s>1$, the initial values $X_0,\ldots,X_{s-1}$ are required.
    If they are not directly available, as usually happens, some approximations can be employed. These must be sufficiently accurate to maintain the desired order of convergence of the original BDF\@. To this end, a classic approach sees the computation of the $(s-j)$-th initial value by a BDF of order $s-j$ for $j=s-1,\ldots,1$.
    Starting from $j=s-1$, we thus compute $X_1$ with a BDF of order 1, and this process is then applied recursively to generate all the required initial values. See~\cite[Algorithm 2]{baranetAl24} for the details of this wind-up procedure.


\section{On the numerical solution of the resulting gAREs}\label{sec:solver4gAREs}
In this section, we revise some standard approaches for the solution of~\eqref{eq:gAREs}. Even though factorization-based methods computing highly accurate solutions exist~\cite{BinietAl2011}, we will employ iterative methods only. This is motivated by a number of factors. First, the problem dimension $n$ could be too large to allocate the dense solution $X_{k}$ in full format. Leveraging the results in, e.g.,~\cite[Section 2]{BennerBujanovic2016}, the exact solution to~\eqref{eq:gAREs} may present a fast decay in its singular values, thus opening the door for accurate low-rank approximations.
In this work, we represent the numerical solution to~\eqref{eq:gAREs} in an $LDL^{\tran}$-format, namely $X_{k}\approx L_{k}D_{k}L_{k}^{\tran}$, $L_{k}\in\mathbb{R}^{n\times r_{k}}$, $D_{k}\in\mathbb{R}^{r_{k}\times r_{k}}$, where we expect $r_{k}$ to be much smaller than $n$. The latter feature makes the storage of $L_{k}$ and $D_{k}$ affordable.
In addition to lowering the memory requirements, working with objects in low-rank format can be beneficial in terms of running time, especially when we need to solve long sequences of equations like in our context. This is often the case also for moderate $n$, when storing full $n\times n$ matrices is indeed possible. Second, as often happens when solving differential problems, the overall error achieved by the computed solution is largely dominated by the discretization error; wasting resources in computing accurate algebraic solutions to~\eqref{eq:gAREs} is thus not necessary. Adopting iterative schemes for~\eqref{eq:gAREs} allows for a trade-off between accuracy of the solution and computational efforts.

\subsection{The Newton-Kleinmann method}\label{sec:NK}

The first method we adopt in our numerical tests is the inexact low-rank Newton-Kleinman (iNK) method~\cite{BenneretAl2016}. Given the gARE~\eqref{eq:gAREs} at time $t_{k}$, iNK solves this equation by a Newton-Kleinman iteration where the related linear step, a Lyapunov equation, is approximated by the low-rank ADI (LR-ADI) method~\cite{BenneretAl2008,BenneretAl2013}. More precisely, given an initial guess $X_{k}^{(0)}$ in factored form $X_{k}^{(0)}= L_{k}^{(0)}D_{k}^{(0)}(L_{k}^{(0)})^{\tran}$ for the Newton-Kleinman sequence, the $\ell$\/th iteration of iNK sees the solution of the following generalized Lyapunov equation
\begin{equation}\label{eq:lyap_iNK}
        \left(G_{k}^{(\ell)}\right)^{\tran}Z_{k}^{(\ell)}\widehat E_{k}+ \widehat E_{k}^{\tran}Z_{k}^{(\ell)}G_{k}^{(\ell)}+P_{k}^{(\ell)}S_{k}^{(\ell)}\left(P_{k}^{(\ell)}\right)^{\tran}=0,
\end{equation}
where
\begin{eqnarray}
\label{eq:lyap_coeff}
    G_{k}^{(\ell)}:=&\widehat{A}_{k}^{(\ell)}-\widehat B_{k}\widehat R_{k}^{-1}\widehat B_{k}^{\tran}X_{k}^{(\ell-1)}\widehat E_{k},\notag\\
    P_{k}^{(\ell)}:=&\left[\widehat{C}_{k}^{\tran},\widehat E_{k}^{\tran}X_{k}^{(\ell-1)}\widehat B_{k}\right],\notag\\
        S_{k}^{(\ell)}:=&\begin{bmatrix}
            \widehat{Q}_{k} & \\
            & \widehat{R}_{k}^{-1}
        \end{bmatrix}.
    \end{eqnarray}
    Equation~\eqref{eq:lyap_iNK} is solved inexactly by LR-ADI which provides a low-rank approximation of the form
\[Z_{k}^{(\ell)}\approx V_{k}^{(\ell)}Y_{k}^{(\ell)}\left(V_{k}^{(\ell)}\right)^{\tran}.\]
At this point, the line search is applied. In particular, a parameter $\lambda_\ell$ is computed and the low-rank factors $L_{k}^{(\ell)}$ and $D_{k}^{(\ell)}$ of the current solution provided by iNK are such that
    \begin{gather}
L_{k}^{(\ell)}D_{k}^{(\ell)}\left(L_{k}^{(\ell)}\right)^{\tran}\approx \notag\\
\left[V_{k}^{(\ell)},L_{k}^{(\ell-1)}\right]
\begin{bmatrix}
    (1-\lambda_\ell) Y_{k}^{(\ell)} & \\
    & D_{k}^{(\ell-1)}
\end{bmatrix}
\left[V_{k}^{(\ell)},L_{k}^{(\ell-1)}\right]^{\tran}.   \label{eq:updatesol_linesearch}
    \end{gather}
See~\cite[Section 3.2]{BenneretAl2016} for more details on the computation of $\lambda_\ell$.
Notice that the matrix on the left-hand side above is never fully allocated, but its low-rank format is maintained and exploited to compute all the quantities required by the following iNK iteration, if needed. Similarly, the low-rank format of $X_{k}^{(\ell-1)}$ can be used when solving shifted linear systems with $G_{k}^{(\ell)}$ within LR-ADI\@; see, e.g.,~\cite{SchS24}.

While the use of LR-ADI for~\eqref{eq:lyap_iNK} allows for an inexact, thus cheaper, solution of the linear step, the computed $L_{k}^{(\ell)}$ and $D_{k}^{(\ell)}$ still need to meet a certain level of accuracy to maintain the convergence of the overall iNK scheme; see~\cite[Section 3.3]{BenneretAl2016}.

The iNK method for~\eqref{eq:gAREs} is reported in Algorithm~\ref{Alg:iNK}, where $\mathcal{R}_{k}(\cdot)$ denotes the generalized Riccati operator in~\eqref{eq:gAREs}; namely, equation~\eqref{eq:gAREs} can be written as
$\mathcal{R}_{k}(X_{k})=0$.
Notice that, thanks to the low rank of $L_{k}^{(\ell-1)}D_{k}^{(\ell-1)}(L_{k}^{(\ell-1)})^{\tran}$, also the residual matrix
$\mathcal{R}_{k}(L_{k}^{(\ell-1)}D_{k}^{(\ell-1)}(L_{k}^{(\ell-1)})^{\tran})$
can be written in a low-rank format. As customary in low-rank methods, this can be exploited to compute its norm, either Frobenius or spectral norm, required by the stopping criterion in line~\ref{iNK_line_stopping} of Algorithm~\ref{Alg:iNK}.

\begin{algorithm}[t]
\DontPrintSemicolon\SetAlgoLined%
\caption{Inexact Newton-Kleinman method for~\eqref{eq:gAREs}\label{Alg:iNK}.}
\KwIn{Coefficient matrices $\widehat A_{k}$, $\widehat E_{k}$, $\widehat B_{k}$, $\widehat R_{k}$,
$\widehat C_{k}$, and $\widehat Q_{k}$, initial guess low-rank factors $L_{k}^{(0)}$, $D_{k}^{(0)}$,
$\mathtt{tol}>0$, max no.\ of iterations $\mathtt{maxit}$.}
\For{$\ell = 1, 2, \ldots, \mathtt{maxit}$}{
  \If{$\|\mathcal{R}_{k}(L_{k}^{(\ell-1)}D_{k}^{(\ell-1)}(L_{k}^{(\ell-1)})^{\tran}) \|\leq\mathtt{tol}$\label{iNK_line_stopping}}{
\Return $L_{k}^{(\ell-1)}$ and $D_{k}^{(\ell-1)}$\;
  }
Set $G_{k}^{(\ell)}$, $P_{k}^{(\ell)}$, and $S_{k}^{(\ell)}$ as in~\eqref{eq:lyap_coeff}\;
  Solve~\eqref{eq:lyap_iNK} by LR-ADI to get $V_{k}^{(\ell)}$, $Y_{k}^{(\ell)}$\;
Compute $\lambda_\ell$\;
Update $L_{k}^{(\ell)}$, $D_{k}^{(\ell)}$ as in~\eqref{eq:updatesol_linesearch}\;
    }
\end{algorithm}

\subsection{RADI}\label{sec:RADI}
The second low-rank method we test is RADI~\cite{BenneretAl2018}. This can be seen as a generalization of LR-ADI to handle the solution of gAREs~\eqref{eq:gAREs}. As LR-ADI, also RADI updates the low-rank factors of the approximate solution by solving a shifted linear system at each iteration. The shift parameters involved in such a step are crucial to achieve a fast convergence in terms of number of iterations. A list of strategies to compute performing RADI (and LR-ADI) shifts can be found in the literature; see, e.g.,~\cite{BenneretAl2018,Kuerschner2019,BenneretAl1415}. Identifying the best shifts for RADI when this is employed as inner method for the solution of gAREs stemming from the BDF discretization of DREs is clearly beyond the scope of this work. In the experiments reported in Section~\ref{sec:numerics}, we will limit ourselves to employing off-the-shelf shifts that often perform well when solving gAREs.

The RADI scheme is recalled in Algorithm~\ref{Alg:RADI}. In this pseudocode, we assume all the shifts $\sigma_\ell$'s to be given from the start. However, an on-the-fly computation of the shifts is also possible. Moreover, effective shifts are often complex. The employment of such shifts necessarily introduces complex arithmetic with a consequent increase in the computational cost of Algorithm~\ref{Alg:RADI}. To alleviate this issue, more sophisticated variants of RADI, minimizing the amount of complex arithmetic, are adopted; see~\cite[Algorithm 2]{BenneretAl2018}. 
As discussed for LR-ADI in Section~\ref{sec:NK}, the low rank of the matrix $K\widehat B_{k}^{\tran}$
can be exploited to reduce the cost of the solution of the linear system in line~\ref{RADI_line_linearsystem} of Algorithm~\ref{Alg:RADI}. To conclude, the norms required in line~\ref{RADI_line_stopping} of Algorithm~\ref{Alg:RADI} can be cheaply computed by exploiting the low rank of the involved matrices.

\begin{algorithm}[t]
\DontPrintSemicolon\SetAlgoLined%
\caption{RADI for~\eqref{eq:gAREs}\label{Alg:RADI}.}
\KwIn{Coefficient matrices $\widehat A_{k}$, $\widehat E_{k}$, $\widehat B_{k}$, $\widehat R_{k}$,
$\widehat C_{k}$, and $\widehat Q_{k}$, initial guess low-rank factors $L_{k}^{(0)}$, $D_{k}^{(0)}$,  $\mathtt{tol}>0$, max no.\ of iterations $\mathtt{maxit}$, shifts $\sigma_\ell$'s.}
$P=\widehat{E}_{k}^{\tran}L_{k}^{(0)}D_{k}^{(0)}(L_{k}^{(0)})^{\tran}\widehat B_{k}$\;

$K=P\widehat{R}_{k}^{-\tran}$\;

$W=[\widehat{A}_{k}^{\tran}L^{(0)}_{k}, \widehat{E}_{k}^{\tran}L^{(0)}_{k}, P, \widehat{C}_{k}^{\tran}]$\;


$S=\texttt{blkdiag}([0, D_{k}^{(0)}; D_{k}^{(0)}, 0], -\widehat{R}_{k}^{-1}, \widehat{Q}_{k})$\;

$Y=(D_{k}^{(0)})^{-1}$\;
\For{$\ell = 1, 2, \ldots, \mathtt{maxit}$}{
  \If{$\|WSW^{*} \|\leq\mathtt{tol}\cdot\|\widehat C_{k}^{\tran}\widehat Q_{k}\widehat C_{k}\|$\label{RADI_line_stopping}}{
\Return $L_{k}^{(\ell-1)}$ and $D_{k}^{(\ell-1)}=Y^{-1}$\;
  }
$V=(\widehat A_{k}^{\tran}-K\widehat B_{k}^{\tran}+\sigma_\ell \widehat E_{k}^{\tran})^{-1}WS$\label{RADI_line_linearsystem}

$L_{k}^{(\ell)}=[L_{k}^{(\ell-1)},\sqrt{-2\text{Re}(\sigma_\ell)}V]$\;
$V_B = V^{*}\widehat{B}_{k}\;$
$\widetilde{Y} = S + V_B\widehat{R}_{k}^{-1}V_B^{*}\;$

$Y=\texttt{blkdiag}(Y,\widetilde Y)$\;
$V_Y = -2\text{Re}(\sigma_{\ell})\widehat{E}_{k}^{\tran}V\widetilde{Y}^{-1}\;$

$W=W+V_Y$\;

$K = K+V_{Y}V_{B}\widehat{R}_{k}^{-1}$\;
    }
\end{algorithm}

\section{Numerical results}\label{sec:numerics}
\setlength{\figwidth}{.9\linewidth}
\setlength{\figheight}{.9\figwidth}
\pgfplotsset{
  every axis/.append style={
    height = \figheight,
    width = \figwidth,
    xminorticks = false,
    yminorticks = true,
    xminorgrids = false,
    yminorgrids = true,
    xmajorgrids = false,
    ymajorgrids = true,
    minor grid style={black!10},
    legend entries = {iNK(0), iNK(1), RADI(0), RADI(1)}
  }
}
The goal of this work is to provide a comprehensive and clear comparison of
state-of-the-art methods for large-scale non-autonomous DREs, with a particular focus on the use of the inner gARE solver. The latter is selected to be either the iNK method with line search or RADI\@.
Both these methods have been tested by using two different initialization strategies: a zero initial guess ($X_{k}^{(0)}=0$ for every $k>0$) or with the solution computed at the previous time step ($X_{k}^{(0)}=X_{k-1}$). In the results that follow, the notation iNK(0) and iNK(1) will be used to denote iNK with zero/non-zero initial guess, respectively. Similarly, we will use the notation RADI(0) and RADI(1) for RADI\@. For any time step $k$ and iNK iteration $\ell$, a zero initial guess has always been adopted when solving~\eqref{eq:lyap_iNK} by LR-ADI. Indeed, in~\cite{SchS24} the authors reported how employing a nonzero initial guess for LR-ADI for Lyapunov equations stemming from a Newton-Kleinman  scheme is not really beneficial, in general, while tackling those for Rosenbrock discretization of DREs benefit greatly. As we expect iNK to suffer the same also in our context, we tested only the case of zero LR-ADI initial guesses.

We mention that we also tested iNK \emph{without} line search. However, its performance was not competitive, neither with a zero initial guess nor with a warm-start, compared to the other schemes. Therefore, iNK with no line search has been excluded from the analysis reported here. However, a comparison involving RADI(0) and iNK without line search in the context of large-scale \emph{autonomous} DREs can be found in~\cite{amslaurea36666}.

In order to provide a comprehensive overview of the performance of the inner solvers, we consider four different datasets, as well as BDF methods of orders ranging from 1 to 4.

\subsection{Software and hardware setups}\label{sec:soft}
The experiments reported here have been executed on a machine with 2 AMD EPYC\textsuperscript{\textsf{TM}} 9554 64-core processors running at 2.67 GHz and equipped with 2 TB total main memory. Using the SLURM scheduling system, our jobs were restricted to 8 cores and 280 GB of main memory. The computer is running on Ubuntu 24.04.4~LTS, and we used MATLAB\textsuperscript{\textcircled{R}} R2025a (version 25.1.0.2943329) for the experiments.

\subsection{Code and data availability}
For three out of the four datasets considered in our experiments, the coefficient matrices are generated by the function \texttt{mess\_get\_linear\_rail} from the \texttt{M.E.S.S.} Toolbox, which now supports both the \texttt{LTI} and \texttt{LTV} options for autonomous and non-autonomous problems, respectively; see Section~\ref{sec:cooling} and~\cite{fenics_rail_mpi} for further details on these datasets.

The fourth dataset consists of synthetic data constructed to satisfy the standard assumptions, ensuring the existence and uniqueness of the problem's solution while at the same time providing a more challenging test case for the considered inner solvers.
Furthermore, for each dataset, we perform simulations for several problem dimensions $n$, while keeping the number of control inputs and outputs fixed at $p=7$ and $q=6$, respectively, unless stated otherwise.
Finally, throughout all non-autonomous experiments, only the coefficient matrices $A$ and $B$ depend on time, as discussed in Sections~\ref{sec:cooling} and~\ref{sec:synt}. Therefore, the mass matrix $E$ remains constant, and equation~\eqref{eq:main} reduces to the case $\dot{E}=0$.

The main routines employed in our experiments are \texttt{mess\_bdf\_dre}, used to solve DREs via BDF discretization, together with either \texttt{mess\_lrnm} or \texttt{mess\_lrradi}, making use of the low-rank structure of the problem via the iNK and RADI methods, respectively. The initialization strategy at each BDF time step, namely zero initialization or warm start, is controlled through the \texttt{opts} structure of the BDF solver.
The same structure is used to specify the stopping criteria for the iterative methods. In particular, the maximum number of iterations is set to $20$ for iNK, $100$ for the ADI iteration arising within iNK, and $150$ for RADI. The remaining stopping criteria are based on the relative residual and the relative difference between two consecutive iterates.
More precisely, the relative residual is defined as the Frobenius norm of the residual of~\eqref{eq:main} evaluated at the current iterate $X_k^{(\ell)}$, normalized by $\|C^{\tran}QC\|_\mathsf{F}$. The relative difference is defined as
\[
\frac{\left|\left|X_k^{(\ell)}-X_k^{(\ell-1)}\right|\right|_\mathsf{F}}
     {\left|\left|X_k^{(\ell)}\right|\right|_\mathsf{F}}.
\]

For both iNK and RADI, the corresponding tolerances are set to $10^{-12}$ and $10^{-16}$, respectively. For the ADI iteration, slightly stricter tolerances are employed, namely $10^{-14}$ for the relative residual and $10^{-16}$ for the relative difference.

As already discussed in Section~\ref{sec:solver4gAREs}, the convergence behavior of the inner solvers strongly depends on the choice of the shift parameters, which are employed indirectly by iNK through the ADI iteration and directly by RADI\@. In our experiments, we tested different shift selection strategies, including Wachspress~\cite{ADIshifts} and Penzl~\cite{Pen00a} shifts. However, these approaches did not provide significant improvements over projection-based shifts~\cite{Kuerschner2019}. Therefore, we only report results obtained with the latter strategy. All routines mentioned above are part of the \texttt{M.E.S.S.} Toolbox~\cite{MMESS}. Our work is based on its version 3.1, and all changes required for this work will be included in the upcoming version 3.2.

For all numerical experiments, we consider the time interval $[t_0,t_f]=[0,720]$, which, in the case of the RAIL examples, corresponds to $7.2\,\mathrm{s}$, since in the simulation model the time variable is scaled by a factor of $10^2$. Moreover, we discretize the time interval using $10$ time steps per unit of time, resulting in a total of $7200$ time steps; hence, $7200$ generalized algebraic Riccati equations must be solved to obtain the full solution trajectory. Finally, throughout all experiments, we assume that the considered problems satisfy~\eqref{eq:main} with terminal condition $X(t_f)=0$.

To quantify the computational effort, we count the number of linear systems of the form
\[
\left(\widetilde{A}+\sigma_\ell \widetilde{E}\right)V=\widetilde{W},
\]
arising within the iterations of the different algorithms, where $\widetilde{A}$, $\widetilde{E}$, and $\widetilde{W}$ denote solver-dependent matrices. This quantity does not necessarily coincide with the number of iterations, since both ADI and RADI may exploit complex shifts, allowing two update steps within a single iteration. Nevertheless, it provides a meaningful measure of computational cost, as the solution of such linear systems constitutes the dominant computational bottleneck.


\subsection{Rail profile}\label{sec:cooling}
The first two datasets we consider arise from an optimal rail profile cooling problem~\cite{10.1007/3-540-27909-1_19}. Briefly, this problem consists of determining an optimal cooling strategy, such as cooling water spray temperatures, locations, or timing, so that the temperature distribution along the rail evolves according to a desired profile while satisfying physical and technological constraints.

The problem is modeled by a system of partial differential equations that, after spatial discretization, gives rise to the following state-space system:
\begin{align} \label{eq:state-space}
E\dot{x} &= \left(\frac{\lambda}{c\rho}P+\frac{\gamma}{c\rho}E_\Gamma\right)x
            + \frac{\gamma}{c\rho}\bar{B}u,
            \qquad t>0,\,x(0)=x_0,\\
y &= Cx,\notag
\end{align}
where $x$ represents the temperature distribution along the rail,
$c = 7620 \frac{m^2}{s^2\,^\circ\mathrm{C}}$ is the specific heat capacity,
$\lambda = 26.4 \frac{kg\,m}{s^3\,^\circ\mathrm{C}}$ is the thermal conductivity,
$\rho = 654 \frac{kg}{m^3}$ is the density of the rail profile,
and $\gamma = 7.0164 \frac{kg}{s^3\,^\circ\mathrm{C}}$.
Moreover, $E$, $P$, $E_\Gamma$, and $\bar{B}$ are matrices specific to the problem.
The optimal control problem is then solved by the Linear Quadratic Regulator (LQR) approach, which consists in the minimization of a quadratic cost function subject to the linear state-space system~\eqref{eq:state-space}.

For our purposes, the relevant observation is that the solution of such an optimization problem depends on the solution $X$ of the DRE~\eqref{eq:main}, where
$A=-(\alpha P+\beta E_\Gamma)$, with $\alpha=\frac{\lambda}{c\rho},\ \beta=\frac{\gamma}{c\rho}$,
and
$B=\beta\bar{B}$.
For more details on the problem background, see~\cite{Eppler2001}. Further details on the specific data generation can be found in~\cite{Saa09}

The first dataset, denoted by \texttt{RAIL\_LTI}, corresponds to the autonomous version of the cooling problem. The second dataset, \texttt{RAIL\_LTV}, considers the corresponding non-autonomous case taken from~\cite{Lan17}, i.e.
\[
\alpha=\alpha(t)=\left(0.75\sin\!\left(2\pi \frac{t}{720}\right)+1\right)\frac{\lambda}{c\rho},
\]
so that $A=A(t)$.
The matrices $E$, $B$, and $C$ are problem-specific and are the same for all two datasets.

In the following, we investigate whether the behavior of the considered solvers changes significantly with respect to both the problem dimension $n$ and the BDF order $s$. We tested the inner solvers for different BDF discretization orders, from BDF1 ($s=1$) to BDF4 ($s=4$). While higher-order discretizations provide more accurate approximations and may improve convergence, especially when warm-start strategies are employed,  they can also lead to handling matrices with larger inner-factor dimensions and, therefore, to larger computational cost.
In our vast numerical tests, BDF2 appeared to provide the best trade-off between accuracy and computational cost. This, along with the sake of simplicity and brevity, made us decide to focus on reporting more detailed information only for the BDF method of order $s=2$.

As for the problem dimension, we considered a range of problem dimensions $n\in\{109,\,371,\,1357,\,5177\}$. However, most of the plots reported here focus on the largest dimension considered for these datasets, namely $n=5177$. Indeed, for smaller dimensions, iNK and RADI tended to exhibit similar performance, making their relative advantages less evident.

\subsubsection{\texttt{RAIL\_LTI}}\label{rail_lti}
Although the main focus of our analysis is in the non-autonomous setting, we begin by presenting results for the autonomous case in order to highlight one of the novelties of this work, namely the use of a warm-start initialization for both iNK and RADI in the solution of DREs, as opposed to the standard zero initial guess commonly adopted in the literature.

The rationale behind this choice is that we are solving a differential equation, and we therefore expect the solution to be sufficiently smooth in time, at least differentiable. In particular, when the time grid is sufficiently fine, consecutive solutions are not expected to differ significantly from each other. This motivates initializing the iterative solvers at each time step with the solution computed at the previous time step so as to supposedly start the iteration closer to the expected solution.

In Table~\ref{tab:RAIL_LTI_5177} we present an overall comparison of the different solvers for various BDF discretization orders and a system of dimension $n=5177$. The table reports the runtimes, the average number of iterations, the average rank of the computed solutions, and the total number of linear systems solved throughout the entire computation.

\begin{figure}[tbp]
    \centering
  \tikzexternalenable%
  \tikzsetnextfilename{runtime_RAIL_LTI}%
  \begin{tikzpicture}
  \begin{axis}[
    ybar,
    bar width=1.25ex,
    xlabel={problem size},
    ylabel={runtime [min]},
    log basis y=10,
    ymode=log,
    legend pos={north west},
    enlarge x limits={0.15, 0.15},
    enlarge y limits={value=0.275, upper},
    symbolic x coords={109, 371, 1357, 5177},
    xtick = data,
    nodes near coords,
    every node near coord/.append style={rotate=90, anchor=west, font=\scriptsize, black},
    visualization depends on={y \as \rawy},
    nodes near coords={\pgfmathparse{10^\rawy}\pgfmathprintnumber{\pgfmathresult}}
    ]

    \addplot[style={mycolor1, fill=mycolor1}] coordinates {(109,2.1063) (371,6.0537) (1357,30.5825) (5177,355.0862)};
    \addlegendentry{iNK(0)}

    \addplot[style={mycolor2, fill=mycolor2}] coordinates {(109,1.3815) (371,4.2987) (1357,22.6027) (5177,294.0668)};
    \addlegendentry{iNK(1)}

    \addplot[style={mycolor3, fill=mycolor3}] coordinates {(109,1.1435) (371,3.3485) (1357,16.6880) (5177,135.5475)};
    \addlegendentry{RADI(0)}

    \addplot[style={mycolor4, fill=mycolor4}] coordinates {(109,1.2605) (371,4.1869) (1357,15.0762) (5177,80.855)};
    \addlegendentry{RADI(1)}

  \end{axis}
\end{tikzpicture}%
  \tikzexternaldisable%

    \caption{Runtime comparison for BDF order $s=2$. Timings are reported in minutes. (\texttt{RAIL\_LTI})}%
    \label{fig:runtime_RAIL_LTI}
\end{figure}
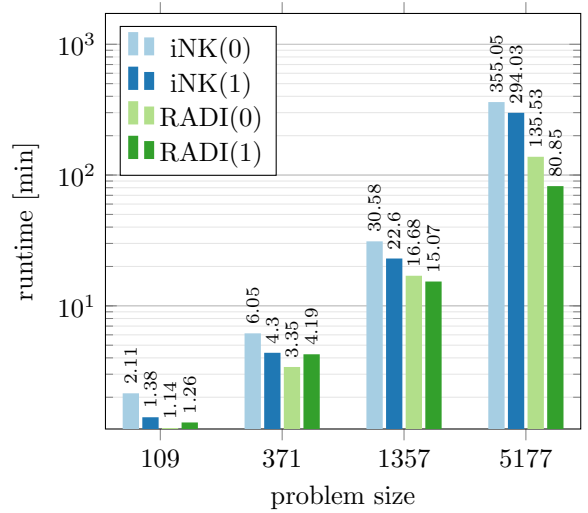

\begin{table}[tbp]
\centering
\caption{Performance comparison for $n=5177$. (\texttt{RAIL\_LTI})}%
\label{tab:RAIL_LTI_5177}

\resizebox{\linewidth}{!}{%
\setlength{\tabcolsep}{4pt}

\begin{tabular}{cccccr}
\toprule
\textbf{BDF} & \textbf{Solver} & \textbf{Runtime} & \textbf{Avg.\ iter.} & \textbf{Avg.\ rank} & \textbf{Lin.\ solves}\\
\midrule

\multirow{4}{*}{$s=1$}
& iNK(0)  & 7h 24 min & 2.0 & 98.50 & 14398 \\
& iNK(1)  & 6h 16 min & 1.1 & 98.50 & 7605 \\
& RADI(0) & 2h 48 min & 8.0 & 98.50 & 57593 \\
& RADI(1) & 1h 20 min & 3.9 & 98.50 & 27893\\
\midrule

\multirow{4}{*}{$s=2$}
& iNK(0)  & 5h 55 min & 2.0 & 97.14 & 14398 \\
& iNK(1)  & 4h 54 min & 1.1 & 97.13 & 7606 \\
& RADI(0) & 2h 16 min & 7.0 & 97.14 & 50395 \\
& RADI(1) & 1h 21 min & 3.5 & 97.14 & 25001 \\
\midrule

\multirow{4}{*}{$s=3$}
& iNK(0)  & 5h 10 min & 2.0 & 96.34 &  14418\\
& iNK(1)  & 4h 16 min & 1.1 & 96.31 & 7626\\
& RADI(0) & 2h 24 min & 7.0 & 96.31 & 50416\\
& RADI(1) & 1h 19 min & 3.0 & 96.31 & 21517 \\
\midrule

\multirow{4}{*}{$s=4$}
& iNK(0)  & 5h 14 min & 2.0 & 95.72 & 14442 \\
& iNK(1)  & 4h 20 min & 1.1 & 95.72 & 7651 \\
& RADI(0) & 2h 6 min & 6.1 & 95.73 & 44302 \\
& RADI(1) & 1h 27 min & 2.8 & 95.72 & 20374 \\
\bottomrule

\end{tabular}
}
\end{table}
We first observe that, across the different BDF orders, the reported quantities exhibit essentially the same qualitative behavior. Second, the warm-start is usually beneficial in terms of average number of iterations; cf. iNK(1) vs iNK(0) and RADI(1) vs RADI(0). In terms of running time, RADI turns out to be more performant than iNK, across all values of $s$ we tested. More remarkably, the total number of linear systems solved decreases as the BDF order increases for RADI. This is likely due to the fact that higher-order discretizations incorporate more information from previously computed solutions, which can be exploited to accelerate convergence. In particular, RADI(1) requires approximately half the number of linear systems needed by RADI(0), and this reduction is consistently reflected in the corresponding runtimes. In case of iNK instead, it looks like the inner ADI is not able to exploit higher BDF orders.
Indeed, the total number of linear systems solved within iNK(0) and iNK(1) is pretty much independent of $s$. Notice that the solution of these many linear systems is the dominant computational cost of iNK. This explains why the number of Newton steps is relatively small while the overall runtime is considerably larger than that of RADI. This is particularly evident in Figure~\ref{fig:runtime_RAIL_LTI}, where we report the runtimes for BDF order $s=2$ and varying the problem dimensions $n$. We observe that iNK(0) is consistently the slowest method. Moreover, as the problem dimension increases, the advantage of iNK(1) over iNK(0) gradually fades away, and the runtimes of the two variants become increasingly similar.
The behavior of the RADI methods is markedly different. For smaller problem dimensions, RADI(0) is slightly faster than RADI(1). However, as $n$ increases, the benefits of the warm-start strategy become more evident, and RADI(1) eventually becomes the fastest method among those considered here.

\begin{figure}[tbp]
    \centering
  \tikzexternalenable%
  \tikzsetnextfilename{RAIL_AUT2_5177_rot}%
  \begin{tikzpicture}
    \pgfplotstableread[col sep=comma]{graphics/data/RAIL2_LTI_5177_last_rot.dat}\mytable

  \begin{semilogyaxis}[%
    ylabel = {solution rank},
    xlabel = {model time},
    legend pos = south east]

    \addplot table[x index = 0, y index = 1] {\mytable};
    \addplot table[x index = 0, y index = 2] {\mytable};
    \addplot table[x index = 0, y index = 3] {\mytable};
    \addplot table[x index = 0, y index = 4] {\mytable};
  \end{semilogyaxis}
\end{tikzpicture}%
  \tikzexternaldisable%

    \caption{Solution ranks over model time, for $s=2$ and $n=5177$. (\texttt{RAIL\_LTI})}%
    \label{fig:rank_RAIL_LTI}
\end{figure}
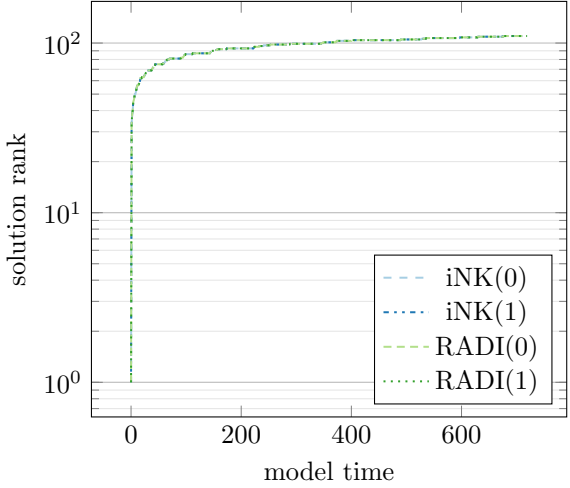

An $s$-dependent trend can be observed for the average rank of the computed solutions. Although all BDF schemes solve the same underlying problem and employ the same truncation tolerance in the column compression procedure, higher-order methods generally produce solutions with lower average rank. This effect is likely due to the fact that lower-order discretizations introduce larger numerical errors, which require additional basis vectors to represent the solution within the prescribed tolerance. Conversely, higher-order BDF schemes are less sensitive to discretization errors and numerical noise and therefore capture the dominant solution components more effectively, resulting in slightly lower-rank approximations.
Furthermore, Figure~\ref{fig:rank_RAIL_LTI} shows the rank of the computed solution as a function of time for $s=2$ and $n=5177$. We observe that, at each time step, all methods produce solutions of essentially identical rank. This is expected, since the same compression tolerances are employed throughout the computations. Differences may arise during the inner iterations used to solve~\eqref{eq:gAREs}, but these disappear after the final compression step, leading to comparable ranks for the converged solutions.

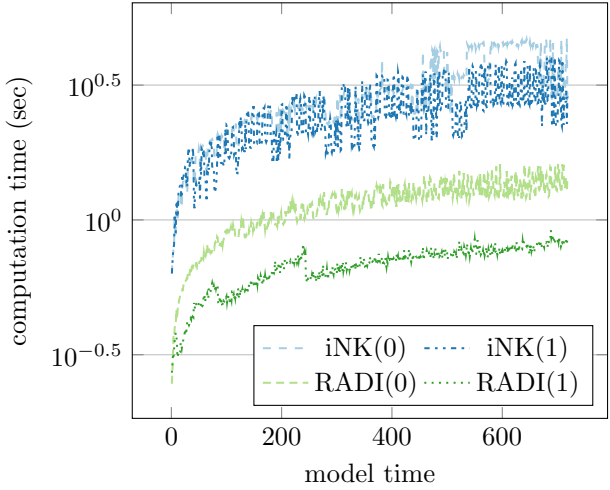
\begin{figure}[tbp]
    \centering
  \tikzexternalenable%
  \tikzsetnextfilename{RAIL_AUT2_5177_tot}%
  \begin{tikzpicture}
  \pgfplotstableread[col sep=comma]{graphics/data/RAIL2_LTI_5177_last_tot.dat}\mytable

  \begin{semilogyaxis}[%
    ylabel = {computation time (sec)},
    xlabel = {model time},
    legend columns=2,
    legend pos = south east
    ]

    \addplot table[x index = 0, y index = 1] {\mytable};
    \addplot table[x index = 0, y index = 2] {\mytable};
    \addplot table[x index = 0, y index = 3] {\mytable};
    \addplot table[x index = 0, y index = 4] {\mytable};
  \end{semilogyaxis}
\end{tikzpicture}%
  \tikzexternaldisable%

    \caption{Inner solver time (in [s]) over model time, for $s=2$ and $n=5177$. (\texttt{RAIL\_LTI})}%
    \label{fig:comptime_RAIL_LTI}
\end{figure}

Finally, Figure~\ref{fig:comptime_RAIL_LTI} illustrates the evolution of the computational time along the simulation horizon for $s=2$ and $n=5177$. All methods exhibit a similar qualitative behavior: the computational cost is relatively low during the first portion of the time interval, increases around the middle of the horizon, and then tends to stabilize in the final part of the simulation. This clearly matches the increment of the rank of the solution. Indeed, working with factors with larger rank, as happens for large $t$ in this example, increases the costs of all the solvers tested here.
The main difference lies in the relative magnitude of the computational costs. The two Newton variants and RADI(1) follow a similar trend, whereas RADI(0) consistently requires less computational time per time step.









\subsubsection{\texttt{RAIL\_LTV}}
The second dataset we consider corresponds to a non-autonomous problem. As described in Section~\ref{sec:cooling}, the \texttt{RAIL\_LTV} dataset uses the same coefficient matrices as \texttt{RAIL\_LTI}, with the only modification being the introduction of a time-dependent coefficient matrix $A$ through the factor $\alpha=\alpha(t)$.

In a time-varying setting, one would naturally expect the problem to become more challenging, since the data in~\eqref{eq:gAREs} change at every BDF step. However, this is not necessarily the case, as the overall difficulty depends on the evolution of the underlying dynamics, which may in some situations even facilitate convergence.

We report the same quantities as in Section~\ref{rail_lti}.
We first discuss Table~\ref{tab:RAIL_LTV_5177} compared with Table~\ref{tab:RAIL_LTI_5177} of the previous section.
\begin{table}[tbp]
\centering
\caption{Performance comparison for $n=5177$. (\texttt{RAIL\_LTV})}%
\label{tab:RAIL_LTV_5177}

\resizebox{\linewidth}{!}{%
\setlength{\tabcolsep}{4pt}

\begin{tabular}{cccccr}
\toprule
\textbf{BDF} & \textbf{Solver} & \textbf{Runtime} & \textbf{Avg.\ iter.} & \textbf{Avg.\ rank} & \textbf{Lin.\ solves}\\
\midrule

\multirow{4}{*}{$s=1$}
& iNK(0)  & 8h 30 min & 2.0 & 100.00 & 14398 \\
& iNK(1)  & 7h 21 min & 1.5 & 100.03 & 10658 \\
& RADI(0) & 3h 3 min & 7.5 & 100.00 & 54218\\
& RADI(1) & 5h 44 min & 6.0 & 100.02 & 43256 \\
\midrule

\multirow{4}{*}{$s=2$}
& iNK(0)  & 6h 25 min & 2.0 & 96.23 & 14398 \\
& iNK(1)  & 5h 42 min & 1.5 & 96.23 & 10660 \\
& RADI(0) & 2h 35 min & 6.8 & 96.27 & 48710 \\
& RADI(1) & 5h 5 min & 5.6 & 96.26 & 40328 \\
\midrule

\multirow{4}{*}{$s=3$}
& iNK(0)  & 4h 58 min & 2.0 & 95.55 &14418\\
& iNK(1)  & 4h 18 min & 1.5 & 95.55 & 10679\\
& RADI(0) & 2h 22 min & 6.5 & 95.48 &  47000 \\
& RADI(1) & 4h 5 min & 5.1 & 95.55 & 36974\\
\midrule

\multirow{4}{*}{$s=4$}
& iNK(0)  & 4h 45 min & 2.0 & 94.95 & 14442 \\
& iNK(1)  & 4h 7 min & 1.5 & 94.95 & 10707 \\
& RADI(0) & 2h 19 min & 6.2 & 94.95 & 44496 \\
& RADI(1) & 3h 45 min & 4.8 & 94.95 & 34502 \\
\bottomrule

\end{tabular}
}
\end{table}

\begin{figure}[tbp]
  \centering
  \tikzexternalenable%
  \tikzsetnextfilename{RAIL2_5177_rot}%
  \begin{tikzpicture}
    \pgfplotstableread[col sep=comma]{graphics/data/RAIL2_LTV_5177_last_rot.dat}\mytable

  \begin{semilogyaxis}[%
    ylabel = {solution rank},
    xlabel = {model time},
    legend pos = south east]

    \addplot table[x index = 0, y index = 1] {\mytable};
    \addplot table[x index = 0, y index = 2] {\mytable};
    \addplot table[x index = 0, y index = 3] {\mytable};
    \addplot table[x index = 0, y index = 4] {\mytable};
  \end{semilogyaxis}
\end{tikzpicture}%
  \tikzexternaldisable%

  \caption{Solution ranks over model time, for $s=2$ and $n=5177$. (\texttt{RAIL\_LTV})}%
  \label{fig:rank_RAIL_LTV}
\end{figure}
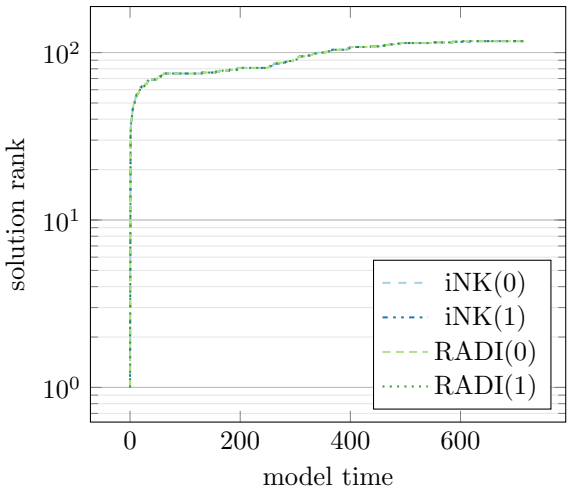

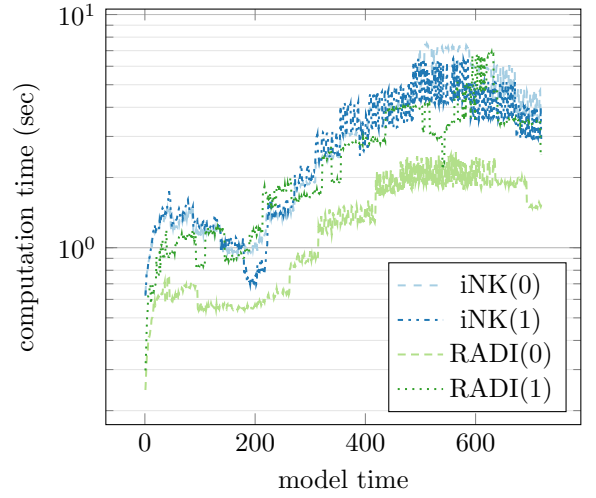
\begin{figure}[tbp]
    \centering
  \tikzexternalenable%
  \tikzsetnextfilename{RAIL2_5177_tot}%
  \begin{tikzpicture}
  \pgfplotstableread[col sep=comma]{graphics/data/RAIL2_LTV_5177_last_tot.dat}\mytable

  \begin{semilogyaxis}[%
    ylabel = {computation time (sec)},
    xlabel = {model time},
    legend pos = south east
    ]

    \addplot table[x index = 0, y index = 1] {\mytable};
    \addplot table[x index = 0, y index = 2] {\mytable};
    \addplot table[x index = 0, y index = 3] {\mytable};
    \addplot table[x index = 0, y index = 4] {\mytable};
  \end{semilogyaxis}
\end{tikzpicture}%
  \tikzexternaldisable%

    \caption{Inner solver time (in [s]) over model time, for $s=2$ and $n=5177$. (\texttt{RAIL\_LTV})}%
    \label{fig:comptime_RAIL_LTV}
\end{figure}

Generally, we observe that all methods tend to require more iterations than in the autonomous case. This behavior is expected. Indeed, the coefficients of~\eqref{eq:gAREs} now depend on time and therefore define a slightly different algebraic problem at each time step. Even though the coefficient matrices evolve smoothly, and the time discretization is sufficiently fine, the problem being solved is no longer exactly the same from one step to the next. As a consequence, information inherited from previous time steps becomes less effective, and additional iterations may be required to achieve convergence. Qualitatively speaking, the performances of iNK(0), iNK(1), and RADI(0) are very similar to those attained for the LTI example. The most significant change concerns RADI(1), whose performance significantly deteriorates in terms of both iteration count and runtime. Although it remains faster than the two iNK variants, RADI(1) becomes approximately two times slower than RADI(0). At the same time, its average number of iterations is still lower than that of RADI(0). This, along with the results reported in Figure~\ref{fig:comptime_RAIL_LTV} for $s=2$ and $n=5177$, suggests that the increased runtime is mainly caused by a higher cost per iteration, likely due to the larger ranks arising during the computation; see also Figure~\ref{fig:rank_RAIL_LTV}.

\begin{figure}[tbp]
  \centering
  \tikzexternalenable%
  \tikzsetnextfilename{runtime_RAIL_LTV}%
  \begin{tikzpicture}
  \begin{axis}[
    ybar,
    bar width=1.25ex,
    xlabel={problem size},
    ylabel={runtime [min]},
    log basis y=10,
    ymode=log,
    legend pos={north west},
    enlarge x limits={0.15, 0.15},
    enlarge y limits={value=0.275, upper},
    symbolic x coords={109, 371, 1357, 5177},
    xtick = data,
    nodes near coords,
    every node near coord/.append style={rotate=90, anchor=west, font=\scriptsize, black},
    visualization depends on={y \as \rawy},
    nodes near coords={\pgfmathparse{10^\rawy}\pgfmathprintnumber{\pgfmathresult}}
    ]

    \addplot[style={mycolor1, fill=mycolor1}] coordinates {(109,2.1483) (371,6.5072) (1357,31.2258) (5177,384.593)};
    \addlegendentry{iNK(0)}

    \addplot[style={mycolor2, fill=mycolor2}] coordinates {(109,1.5157) (371,4.7247) (1357,25.7168) (5177,342.2807)};
    \addlegendentry{iNK(1)}

    \addplot[style={mycolor3, fill=mycolor3}] coordinates {(109,1.1567) (371,3.6367) (1357,16.4262) (5177,154.9265)};
    \addlegendentry{RADI(0)}

    \addplot[style={mycolor4, fill=mycolor4}] coordinates {(109,1.6190) (371,5.7562) (1357,33.6833) (5177,304.9148)};
    \addlegendentry{RADI(1)}

  \end{axis}
\end{tikzpicture}%
  \tikzexternaldisable%

  \caption{Runtime comparison for BDF order $s=2$. (\texttt{RAIL\_LTV}) Timings are reported in minutes.}%
  \label{fig:runtime_RAIL_LTV}
\end{figure}
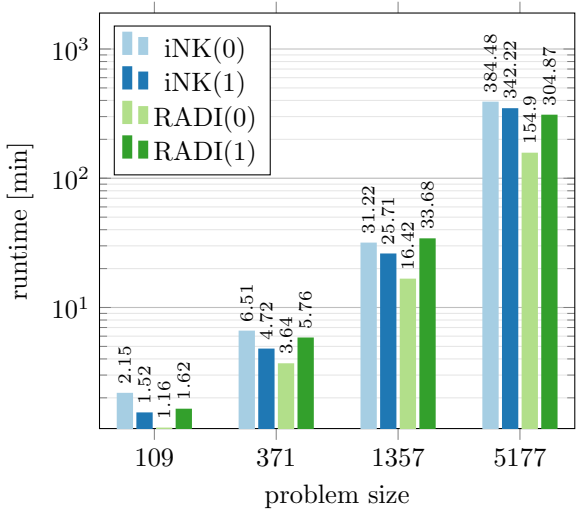

\begin{figure}[tbp]
    \centering
  \tikzexternalenable%
  \tikzsetnextfilename{runtime_RAILR_LTV}%
  \begin{tikzpicture}
  \begin{axis}[
    ybar,
    bar width=1.25ex,
    xlabel={problem size},
    ylabel={runtime [min]},
    log basis y=10,
    ymode=log,
    legend pos={north west},
    enlarge x limits={0.15, 0.15},
    enlarge y limits={value=0.275, upper},
    symbolic x coords={109, 371, 1357, 5177},
    xtick = data,
    nodes near coords,
    every node near coord/.append style={rotate=90, anchor=west, font=\scriptsize, black},
    visualization depends on={y \as \rawy},
    nodes near coords={\pgfmathparse{10^\rawy}\pgfmathprintnumber{\pgfmathresult}}
    ]

    \addplot[style={mycolor1, fill=mycolor1}] coordinates {(109,3.7495) (371,10.1658) (1357,45.3862) (5177,419.4203)};
    \addlegendentry{iNK(0)}

    \addplot[style={mycolor2, fill=mycolor2}] coordinates {(109,2.2283) (371,6.3222) (1357,29.3033) (5177,298.9018)};
    \addlegendentry{iNK(1)}

    \addplot[style={mycolor3, fill=mycolor3}] coordinates {(109,1.2413) (371,3.9618) (1357,15.3823) (5177,131.6745)};
    \addlegendentry{RADI(0)}

    \addplot[style={mycolor4, fill=mycolor4}] coordinates {(109,1.8093) (371,5.7758) (1357,36.5800) (5177,285.6878)};
    \addlegendentry{RADI(1)}

  \end{axis}
\end{tikzpicture}%
  \tikzexternaldisable%

    \caption{Runtime comparison for BDF order $s=2$. (\texttt{RAILR\_LTV}) Timings are reported in minutes.}%
    \label{fig:runtime_RAILR_LTV}
\end{figure}
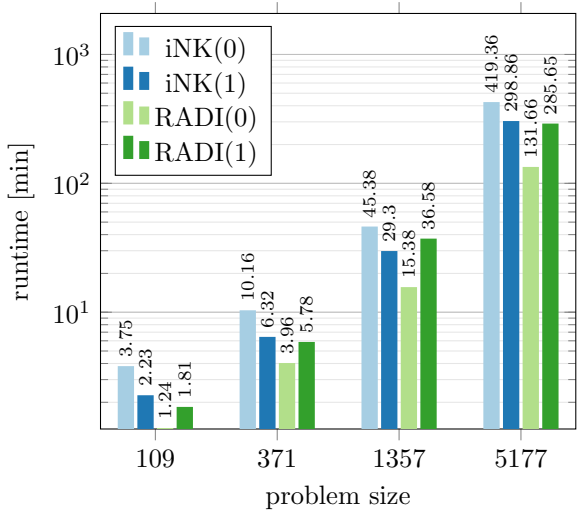

In order to further assess the consistency of the results, we consider an additional dataset, denoted by \texttt{RAILR\_LTV}. This dataset is based on the same data as \texttt{RAIL\_LTV}; however, while in the first two datasets the weighting matrices $Q$ and $R$ are chosen as identity matrices, here we set $R = 10^{-6} I$. This choice brings the problem closer to a nearly singular setting and therefore makes it more challenging to solve.
From a control theory perspective, such a choice corresponds to a strong reduction in the control penalty, meaning that the use of the control input becomes significantly less expensive compared to the state cost. Since the optimal feedback gain depends on $R^{-1}$, this leads to a feedback gain with a substantially larger norm and, consequently, a more aggressive control law. In spite of this difference in the underlying optimal control problem, the numerical results do not reveal any significant new behavior, as they are very similar to those obtained for the \texttt{RAIL\_LTV} dataset. For this reason, we restrict the presentation to the runtime plot in Figure~\ref{fig:runtime_RAILR_LTV}, where we compare the different solvers for BDF order $s=2$.
The only noticeable effect is that iNK(0) appears to be the most sensitive method with respect to this modification of $R$, as reflected in a larger number of Newton steps. This behavior explains its increased computational cost relative to the other methods.

\subsection{Synthetic data}\label{sec:synt}
The last dataset considered in our experiments, denoted by \texttt{SYNT}, corresponds to a non-autonomous problem that does not originate from a specific real-world application. Instead, it is synthetically constructed to provide a challenging yet well-controlled benchmark for the considered numerical methods. The DRE has again the form~\eqref{eq:main}, with terminal condition $X(t_f)=0$. In contrast to the cooling problem discussed in Section~\ref{sec:cooling}, the coefficient matrices are dense. More precisely, we consider
\begin{align*}
A(t) &= \texttt{randn}(n) - \left(\sqrt{n} + \exp(t)\right)I, \\
E &= I, \\
B(t) &= \left(t+1\right)I, \\
C &= \texttt{ones}(1,n),
\end{align*}
with $Q=R=I$, where $n$ denotes the problem dimension, the random matrices are generated using \texttt{rng(1)}, and $I$ denotes the identity matrix of appropriate size.

From the state-space interpretation, this choice of data defines a dynamical system whose stability increases with time, due to the progressively larger negative shift in the matrix $A(t)$. At the same time, the control effectiveness increases with time through the factor $(t+1)$ in $B(t)$. The choice of $C=\texttt{ones}(1,n)$ implies that the output corresponds to the aggregate contribution of all state variables, rather than to a selection of specific components.

Unlike the rail-cooling datasets, this problem does not exhibit any particular sparsity structure and involves dense coefficient matrices throughout the computation. As a consequence, the computational burden is significantly higher, and the methods already encounter difficulties for moderate problem dimensions. For this reason, we restrict our comparison to $n\in\{109,371,1357\}$.
Nevertheless, we expect the qualitative behavior observed in the following experiments to remain representative of larger-scale instances of the same problem class.

\subsubsection{\texttt{SYNT}}
For this synthetic dataset, we observed a rather unexpected behavior when considering the largest problem dimension, $n=1357$, together with a BDF discretization of order $s=4$ and a fine time grid consisting of $\ell=7200$ time steps. Since, in this setting, some considered methods fail to converge, we base our main comparison on the case $\ell=3600$, corresponding to $5$ time discretization points per unit of time, where all methods successfully complete the simulation.

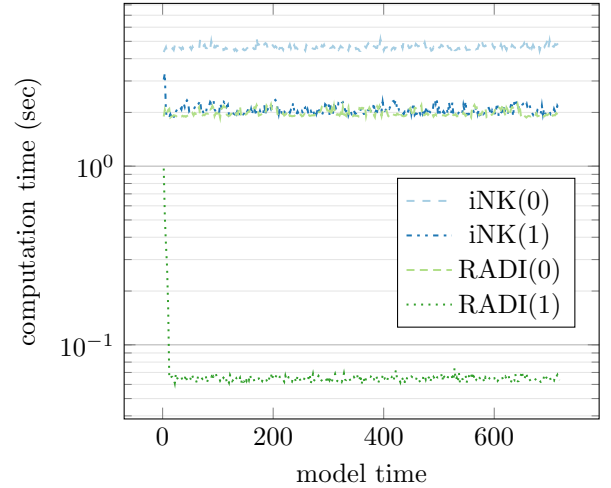
\begin{figure}[tbp]
    \centering
  \tikzexternalenable%
  \tikzsetnextfilename{DATASYNT2_nt5_1357_tot}%
  \begin{tikzpicture}
  \pgfplotstableread[col sep=comma]{graphics/data/DATA2_nt5_1357_last_tot.dat}\mytable

  \begin{semilogyaxis}[%
    ylabel = {computation time (sec)},
    xlabel = {model time},
    legend style={at={(0.95,0.4)},anchor=east}
    ]

    \addplot table[x index = 0, y index = 1] {\mytable};
    \addplot table[x index = 0, y index = 2] {\mytable};
    \addplot table[x index = 0, y index = 3] {\mytable};
    \addplot table[x index = 0, y index = 4] {\mytable};
  \end{semilogyaxis}
\end{tikzpicture}%
  \tikzexternaldisable%

    \caption{Inner solver time (in [s]) over model time, for $s=2$, \(\ell=3600\), and $n=1357$. (\texttt{SYNT})}%
    \label{fig:comptime_SYNT5}
\end{figure}

The first observation arising from Table~\ref{tab:SYNT_5_1357} is the remarkable savings coming from the warm-start strategy. Indeed, the gains achieved by warm-starting are significantly larger than those observed for the previous datasets and are consistent across all considered BDF orders. Although \texttt{SYNT} is a non-autonomous problem, all coefficient matrices depend smoothly on time. Consequently, when a sufficiently fine time grid is employed, the algebraic Riccati equation to be solved at a given time step differs only slightly from the one at the previous step. In this situation, the solution computed at the previous time step provides an extremely accurate initial guess, making warm-starting particularly effective.
\begin{table}[tbp]
\centering
\caption{Performance comparison for $n=1357$ and $\ell=3600$. (\texttt{SYNT})}%
\label{tab:SYNT_5_1357}

\resizebox{\linewidth}{!}{%
\setlength{\tabcolsep}{4pt}

\begin{tabular}{cccccr}
\toprule
\textbf{BDF} & \textbf{Solver} & \textbf{Runtime} & \textbf{Avg.\ iter.} & \textbf{Avg.\ rank} & \textbf{Lin.\ solves}\\
\midrule

\multirow{4}{*}{$s=1$}
& iNK(0) & 5h 28 min & 3.0 & 40.98 & 5398.5 \\
& iNK(1) & 2h 30 min & 1.0 & 40.98 & 1812.5 \\
& RADI(0) & 2h 16 min & 41.0 & 40.98 & 147521 \\
& RADI(1) & 4 min & 1.1 & 40.98 &  4076 \\
\midrule

\multirow{4}{*}{$s=2$}
& iNK(0) & 4h 39 min & 3.0 & 40.97 & 5398.5 \\
& iNK(1) & 2h 5 min & 1.0 & 40.97 & 1814 \\
& RADI(0) & 1h 59 min & 37.0 & 40.96 & 133138 \\
& RADI(1) & 4 min & 1.1 & 40.97 & 4123 \\
\midrule

\multirow{4}{*}{$s=3$}
& iNK(0) & 4h 28 min & 3.0 & 40.87 & 5411.5 \\
& iNK(1) & 2h 1 min & 1.0 & 40.87 & 1837 \\
& RADI(0) & 1h 47 min & 33.9 & 40.86 & 122415 \\
& RADI(1) & 5 min & 1.3 & 40.87 & 4523 \\
\midrule

\multirow{4}{*}{$s=4$}
& iNK(0) & 4h 14 min & 3.0 &  41.19 & 5427 \\
& iNK(1) & 1h 58 min & 1.1 & 41.19 & 1969.5 \\
& RADI(0) & 1h 51 min & 34.8 & 41.19 & 125837 \\
& RADI(1) & 11 min & 3.1 & 41.19 & 11313 \\
\bottomrule

\end{tabular}
}
\end{table}
This behavior is clearly reflected in the results. For iNK, warm-starting yields a speed-up factor of approximately $2.3\times$ together with a reduction by a factor of about $3$ in the number of solved linear systems, corresponding to an overall runtime reduction of roughly $57\%$. The effect is even more pronounced for RADI. In this case, warm-starting achieves a speed-up factor of approximately $29.8\times$ and reduces the number of solved linear systems by a factor of about $32.3$, resulting in a runtime reduction of nearly $97\%$. Consequently, the warm-started version of RADI substantially outperforms all the other methods considered in our experiments. Interestingly, despite the improvements obtained by warm-starting, the performance of iNK(1) remains comparable to that of RADI(0). This suggests that, for problems with very smooth dynamics, the simplicity and efficiency of the RADI iteration can compensate for the additional flexibility provided by the Newton framework.
\begin{figure}[tbp]
    \centering
  \tikzexternalenable%
  \tikzsetnextfilename{runtime_SYNT5}%
  \begin{tikzpicture}
  \begin{axis}[
    ybar,
    bar width=1.25ex,
    xlabel={problem size},
    ylabel={runtime [min]},
    log basis y=10,
    ymode=log,
    legend pos={north west},
    enlarge x limits={0.15, 0.15},
    enlarge y limits={value=0.275, upper},
    symbolic x coords={109, 371, 1357},
    xtick = data,
    nodes near coords,
    every node near coord/.append style={rotate=90, anchor=west, font=\scriptsize, black},
    visualization depends on={y \as \rawy},
    nodes near coords={\pgfmathparse{10^\rawy}\pgfmathprintnumber{\pgfmathresult}}
    ]

    \addplot[style={mycolor1, fill=mycolor1}] coordinates {(109, 4.789) (371, 20.1478) (1357, 278.8182)};
    \addlegendentry{iNK(0)}

    \addplot[style={mycolor2, fill=mycolor2}] coordinates {(109, 1.474) (371, 6.4682 ) (1357, 124.5353)};
    \addlegendentry{iNK(1)}

    \addplot[style={mycolor3, fill=mycolor3}] coordinates {(109, 1.6565) (371, 6.5668) (1357, 119.2922)};
    \addlegendentry{RADI(0)}

    \addplot[style={mycolor4, fill=mycolor4}] coordinates {(109, 0.4968) (371, 1.2017) (1357, 4.2925)};
    \addlegendentry{RADI(1)}

  \end{axis}
\end{tikzpicture}%
  \tikzexternaldisable%

    \caption{Runtime comparison for BDF order $s=2$ and $\ell=3600$. (\texttt{SYNT}) Timings are reported in minutes.}%
    \label{fig:runtime_SYNT5}
\end{figure}
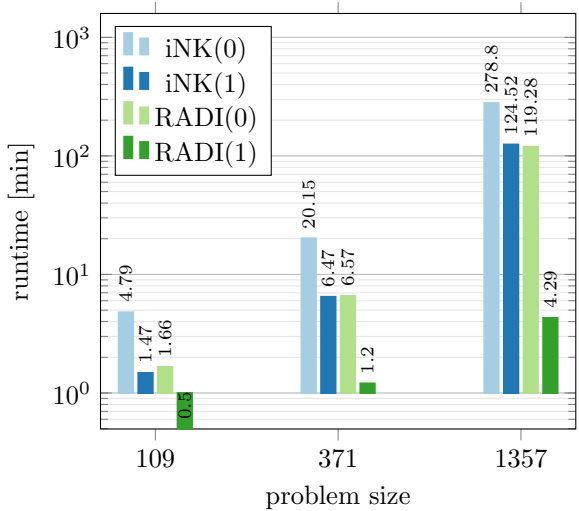
These observations are further supported by Figures~\ref{fig:runtime_SYNT5} and~\ref{fig:comptime_SYNT5}, which display the overall runtimes and the computational time per time step throughout the simulation horizon, respectively. The superior efficiency of RADI(1) is clearly evident, whereas iNK(1) and RADI(0) exhibit comparable performance. As expected, iNK(0) remains the most expensive method over the entire computation.

A second noteworthy aspect concerns the overall computational difficulty of the \texttt{SYNT} dataset. Even for $n=1357$, the runtimes of most methods are comparable to, or even larger than, those observed for the RAIL datasets with $n=5177$. This behavior can be explained by the structure of the coefficient matrices. In the rail-cooling examples, the coefficient matrix $A$ exhibits a high degree of sparsity, which can be exploited to reduce both memory requirements and computational costs. In contrast, in the \texttt{SYNT} dataset, $A$ is dense, leading to more expensive linear algebra operations and a larger computational burden. Remarkably, RADI(1) is affected only marginally by this increased complexity, since the warm-start initialization places the iteration sufficiently close to the solution, requiring only a few additional updates.

Another interesting feature of this dataset is the influence of the BDF order. For all methods except RADI(1), the computational performance improves significantly as the BDF order increases. This behavior is consistent with the observations made for the previous datasets and can be attributed to the larger amount of information inherited from previously computed solutions. In contrast, the performance of RADI(1) remains essentially unchanged, suggesting that the benefits provided by the increasingly accurate initial guess dominate the effects associated with the higher-rank quantities arising in higher-order BDF schemes.

Finally, we discuss the anomalous behavior observed for $\ell=7200$. In this case, none of the methods successfully solve the problem when combined with a BDF discretization of order $s=4$: the Newton schemes fail to converge despite the use of line search, while RADI reaches the prescribed maximum number of iterations. This phenomenon is observed only for BDF4 and only for the finest time grid.
A possible explanation is related to cancellation effects arising in the BDF4 discretization. Since the solution at the current time step is expressed as a linear combination of several previous solutions, difficulties may arise when the time grid is very fine and the exact solution varies only marginally from one step to the next. In this regime, the subtraction of nearly identical quantities may amplify numerical errors and degrade the quality of the resulting algebraic problem.
This interpretation is supported by two observations. First, throughout all previous experiments, the average solution rank decreases as the BDF order increases. In contrast, Table~\ref{tab:SYNT_5_1357} shows a different trend: while the average rank decreases from $s=1$ to $s=3$, it increases again for $s=4$, even exceeding the value observed for $s=1$. This suggests that the inner solvers require higher-rank approximations to capture the relevant features of the solution. Second, the convergence difficulties disappear when the time grid is coarsened from $\ell=7200$ to $\ell=3600$, supporting the hypothesis that the issue is related to the interaction between a very fine temporal discretization and the high-order BDF scheme.
\begin{figure}[tbp]
    \centering
  \tikzexternalenable%
  \tikzsetnextfilename{runtime_SYNT10}%
  \begin{tikzpicture}
  \begin{axis}[
    ybar,
    bar width=1.25ex,
    xlabel={problem size},
    ylabel={runtime [min]},
    log basis y=10,
    ymode=log,
    legend pos={north west},
    enlarge x limits={0.15, 0.15},
    enlarge y limits={value=0.275, upper},
    symbolic x coords={109, 371, 1357},
    xtick = data,
    nodes near coords,
    every node near coord/.append style={rotate=90, anchor=west, font=\scriptsize, black},
    visualization depends on={y \as \rawy},
    nodes near coords={\pgfmathparse{10^\rawy}\pgfmathprintnumber{\pgfmathresult}}
    ]

    \addplot[style={mycolor1, fill=mycolor1}] coordinates {(109, 5.2772) (371, 34.747) (1357, 464.6057)};
    \addlegendentry{iNK(0)}

    \addplot[style={mycolor2, fill=mycolor2}] coordinates {(109, 2.2108) (371, 12.6882) (1357, 198.2608)};
    \addlegendentry{iNK(1)}

    \addplot[style={mycolor3, fill=mycolor3}] coordinates {(109, 2.655) (371, 13.798) (1357, 189.1465)};
    \addlegendentry{RADI(0)}

    \addplot[style={mycolor4, fill=mycolor4}] coordinates {(109, 0.9232) (371, 2.5382) (1357, 8.239)};
    \addlegendentry{RADI(1)}

  \end{axis}
\end{tikzpicture}%
  \tikzexternaldisable%

    \caption{Runtime comparison for BDF order $s=2$ and $\ell=7200$. (\texttt{SYNT})}%
    \label{fig:runtime_SYNT10}
\end{figure}
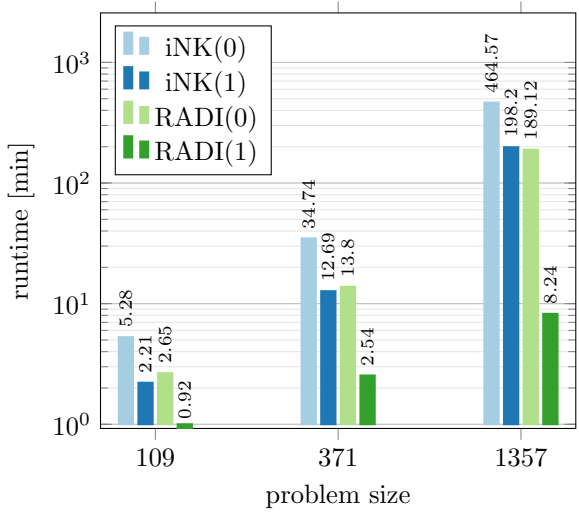
For completeness, the results obtained with $\ell=7200$ are reported in Table~\ref{tab:SYNT_10_1357}, while Figure~\ref{fig:runtime_SYNT10} compares the runtimes obtained for BDF order $s=2$ on the finer grid.

\begin{table}[tbp]
\centering
\caption{Performance comparison for $n=1357$ and $\ell=7200$. (\texttt{SYNT})}%
\label{tab:SYNT_10_1357}

\resizebox{\linewidth}{!}{%
\setlength{\tabcolsep}{4pt}

\begin{tabular}{cccccr}
\toprule
\textbf{BDF} & \textbf{Solver} & \textbf{Runtime} & \textbf{Avg.\ iter.} & \textbf{Avg.\ rank} & \textbf{Lin.\ solves} \\
\midrule

\multirow{4}{*}{$s=1$}
& iNK(0) & 8h 29 min & 3.0 & 40.96 & 10798.5 \\
& iNK(1) & 3h 40 min & 1.0 & 40.96 & 3621 \\
& RADI(0) & 3h 38 min & 35.0 & 40.96 & 251854 \\
& RADI(1) & 8 min & 1.1 & 40.96 & 7825 \\
\midrule

\multirow{4}{*}{$s=2$}
& iNK(0) & 7h 44 min & 3.0 & 40.95 & 10798.5 \\
& iNK(1) & 3h 18 min & 1.0 & 40.95 & 3623 \\
& RADI(0) & 3h 9 min & 29.0 & 40.95 & 208699 \\
& RADI(1) & 8 min & 1.1 & 40.95 &  7800 \\
\midrule

\multirow{4}{*}{$s=3$}
& iNK(0) & 6h 38 min & 3.0 & 40.90 & 10811 \\
& iNK(1) & 2h 45 min & 1.0 & 40.90 & 3649.5 \\
& RADI(0) & 2h 39 min & 27.0 & 40.89 &  194361 \\
& RADI(1) & 9 min & 1.2 & 40.90 &  8373 \\
\bottomrule

\end{tabular}
}
\end{table}

\begin{figure}[tbp]
  \centering
  \tikzexternalenable%
  \tikzsetnextfilename{DATASYNT2_nt5_1357_rot}%
  \begin{tikzpicture}
    \pgfplotstableread[col sep=comma]{graphics/data/DATA2_nt5_1357_last_rot.dat}\mytable

  \begin{semilogyaxis}[%
    ylabel = {solution rank},
    xlabel = {model time},
    legend pos = south east]

    \addplot table[x index = 0, y index = 1] {\mytable};
    \addplot table[x index = 0, y index = 2] {\mytable};
    \addplot table[x index = 0, y index = 3] {\mytable};
    \addplot table[x index = 0, y index = 4] {\mytable};
  \end{semilogyaxis}
\end{tikzpicture}%
  \tikzexternaldisable%

  \caption{Solution ranks over model time, for $s=2$, \(nt=3600\), and $n=1357$. (\texttt{SYNT})}%
  \label{fig:rank_SYNT5}
\end{figure}
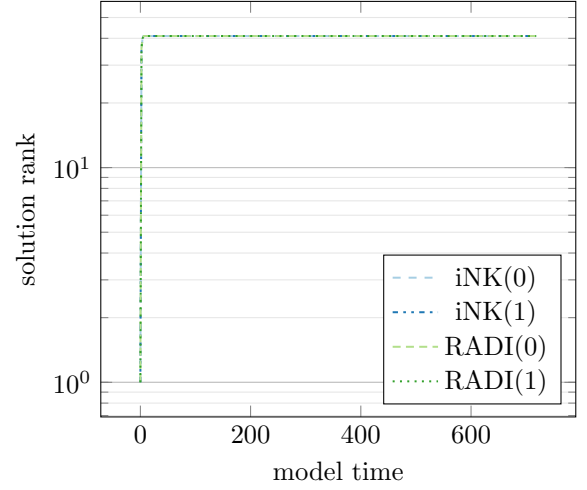

Finally, comparing Figure~\ref{fig:rank_SYNT5} with Figure~\ref{fig:rank_RAIL_LTI}, we observe a markedly different evolution of the solution ranks. In the \texttt{SYNT} dataset, the ranks quickly reach a plateau after only a few time steps, whereas in the \texttt{RAIL\_LTI} example they gradually increase throughout the simulation. This difference reflects the distinct structure of the underlying control problems. 

\section{Conclusions}\label{sec:concl}

In this paper, we investigated the numerical solution of large-scale non-autonomous differential Riccati equations discretized in time by Backward Differentiation Formula (BDF) schemes. In particular, we compared two state-of-the-art low-rank solvers for the generalized algebraic Riccati equations arising at each time step, namely the inexact Newton--Kleinman (iNK) method and RADI, combined with two different initialization strategies: a zero initial guess and a warm-start approach based on the solution computed at the previous time step.

The numerical experiments conducted on several datasets and for different BDF orders reveal a number of consistent trends. First, iNK initialized with a zero matrix is systematically the least competitive approach, both in terms of runtime and computational effort. Although warm-starting significantly improves the performance of iNK by reducing the number of Newton iterations, the resulting gains strongly depend on the smoothness of the underlying problem. In particular, the benefits are substantial for the synthetic dataset, where the coefficient matrices vary smoothly in time, while they become less pronounced for the rail-cooling examples as the problem dimension increases.

Regarding the RADI schemes, they consistently outperform the Newton-based methods across all considered datasets and discretization orders. However, the effectiveness of warm-starting is problem dependent. While it provides substantial improvements for the \texttt{RAIL\_LTI} dataset and exceptionally large gains for the \texttt{SYNT} example, it may also lead to a deterioration of performance, as observed for the \texttt{RAIL\_LTV} problem. This indicates that, in the non-autonomous setting, the quality of a warm-start depends not only on the smoothness of the solution but also on the evolution of the coefficient matrices defining the problem.

Concerning the time discretization, higher-order BDF schemes generally reduce the number of linear systems that must be solved and therefore lead to lower runtimes. On the other hand, they require larger memory resources and may become more sensitive to numerical effects. In particular, for the synthetic dataset we observed convergence difficulties for the fourth-order scheme when a very fine temporal grid was employed. These observations suggest that second-order BDF discretizations provide the most robust compromise between efficiency, memory requirements, and numerical stability.

Overall, the experiments indicate that RADI combined with an appropriate initialization strategy represents the most effective approach among those considered for the numerical solution of large-scale differential Riccati equations. At the same time, the results highlight the importance of carefully balancing the choice of the time discretization scheme and the initialization strategy, especially in the non-autonomous setting.

\section*{Acknowledgments}
The first two authors are members of the Italian INdAM Research group GNCS\@.

\bibliographystyle{siamplain}
\bibliography{ref}

\end{document}